\documentclass[11pt]{article}
\usepackage{amsmath,amssymb,amsthm}

\newtheorem{theorem}{Theorem}[section]
\newtheorem{lemma}[theorem]{Lemma}
\newtheorem{proposition}[theorem]{Proposition}
\newtheorem{corollary}[theorem]{Corollary}

\theoremstyle{plain}

\theoremstyle{definition}
\newtheorem{definition}[theorem]{Definition}

\numberwithin{equation}{section}

\renewcommand{\labelenumi}{\textup{(\theenumi)}}
\begin{document}
\title{On the Markov-Dyck shifts of vertex type}
%\author{Wolfgang Krieger}
%\address{Institute for Applied  Mathematics, University of Heidelberg,
%Im Neuenheimer Feld 294, 69120 Heidelberg, Germany}
%\email{krieger{@@}math.uni-heidelberg.de}
\author{Kengo Matsumoto\\
Department of Mathematics, \\
Joetsu University of Education,\\
Joetsu 943-8512 Japan}
\date{}
%\email{kengo{@@}juen.ac.jp}
\maketitle
\begin{abstract}
For a given finite directed graph $G$, 
there are two types of Markov-Dyck shifts,
the Markov-Dyck shift $D_G^V$ of vertex type
and  
the Markov-Dyck shift $D_G^E$ of edge type.
It is shown that, if $G$ does not have multi-edges,
the former is a finite-to-one factor of the latter,
and 
 they have the same topological entropy.
An expression for the zeta function of a Markov-Dyck shift of vertex type 
is given. 
It is different from that of 
the  Markov-Dyck shift of edge type.
\end{abstract}

%%%%%%%%%%%%%%%%%%%%%%%%%%%%%%%%%%%%%%%%%%%%%%%%%%%%%%                          
\def\C{{{\cal C}}}
\def\V{{{\cal V}}}
\def\det{{{\operatorname{det}}}}
\def\card{{{\operatorname{card}}}}
\def\diag{{{\operatorname{diag}}}}
\def\trace{{{\operatorname{trace}}}}
\def\card{{{\operatorname{card}}}}
\def\N{{ {\mathbb{N}} }}
\def\Z{{ {\mathbb{Z}} }}
\def\OA{{ {\cal O}_A }}
%%%%%%%%%%%%%%%%%%%%%%%%%%%%%%%%%%%%%%% 



\bigskip

%%%%%%%%%%%%%%%%%%%%%%%%%%%%%%%%%%%%
\section{Introduction}
%%%%%%%%%%%%%%%%%%%%%%%%%%%%%%%%%

Let $\Sigma$ be a finite alphabet, and let $\sigma$ be the left shift on
$\Sigma^{\Z}$ defined by
$
 \sigma((x_n)_{ n \in {\Z}})   = (x_{n+1})_{ n \in {\Z}},
$   
$
(x_n)_{ n \in {\Z}} \in \Sigma^{\Z}.
$
For a closed subset $\Lambda \subset \Sigma^\Z$ satisfying 
$\sigma(\Lambda) = \Lambda$,
the topological dynamical system
$(\Lambda, \sigma)$ is called a subshift.
Denote by $B_n(\Lambda)$ the set of all admissible words 
appearing in $\Lambda$ with length $n$,
and by $P_n(\Lambda)$ 
the set of all $n$-periodic points of $(\Lambda,\sigma)$,
respectively.
Then the topological entropy $h_{top}(\Lambda)$
and the zeta function 
$\zeta_\Lambda(z)$ for $(\Lambda,\sigma)$
is defined by
\begin{align}
h_{top}(\Lambda) & 
= \lim_{n \to \infty} \frac{1}{n} \log |B_n(\Lambda)|, \label{eq:entropy}\\
\zeta_\Lambda(z)&
 = \exp(\sum_{n=1}^\infty \frac{|P_n(\Lambda)| z^n}{n}).\label{eq:zeta}
\end{align}
They are crucial topological conjugacy invariants of $(\Lambda,\sigma)$.
 For an introduction to their theory, which belongs to symbolic
dynamics, we refer to \cite{Ki} and \cite{LM}.

W. Krieger in \cite{Kr1} has introduced the Dyck shifts 
from automata theory and language theory in computer science.
They are  non-sofic subshifts defined by Dyck languages.
In 
\cite{HIK,Kr1,Kr2,KMMunster,MathScand2011},
a class of non-sofic subshifts called Markov-Dyck shifts have been studied
(cf. \cite{HK}).
The subshifts are generalization of Dyck shifts 
by using finite directed graphs.
They have recently come to be studied by computer scientists 
(cf. \cite{BBD,CS}).
For a given finite directed graph $G=(V,E)$, 
there are two types of Markov-Dyck shifts,
the Markov-Dyck shift $D_G^V$ of vertex type
and  
the Markov-Dyck shift $D_G^E$ of edge type.
Both of them are not sofic subshifts 
if $G$ is irreducible and not permutive.  
In the papers \cite{HIK,Kr1,Kr2,KMMunster}, 
the Markov-Dyck shifts  mean the Markov-Dyck shifts of edge type.
In  \cite{KMMunster}, formulae of topological entropy and zeta functions for
Markov-Dyck shifts of edge type have been presented.

In the first part of the paper,
we will study relationship between the two types of Markov-Dyck shifts
for finite directed graphs,
the Markov-Dyck shift $D_G^V$ of vertex type
and  
the Markov-Dyck shift $D_G^E$ of edge type.
We will show that, if $G$ does not have multi-edges,
there exists a finite-to-one factor code from
$D_G^E$ to $D_G^V$ (Proposition \ref{prop:factor}). 
The factor code can never yield a topological conjugacy unless the transition matrix of the graph is permutation. 
They have the same topological entropy (Theorem \ref{thm:entropy}).

In the second part of the paper,
we will present a formula of  the zeta function 
of a Markov-Dyck shift of vertex type (Theorem \ref{thm:main}).
The formula is regarded as a generalization of the formula for 
Markov-Dyck shifts of edge type \cite[Theorem 2.3]{KMMunster}.
In the final section,
the zeta function of the  Fibonacci-Dyck shift
of  vertex type will be  presented.
It is different from that of 
the  Fibonacci-Dyck shift
of edge type.
Hence the  Fibonacci-Dyck shift
of  vertex type is not topologically conjugate to
the  Fibonacci-Dyck shift of edge type.

%%%%%%%%%%%%%%%%%%%%%%%%%%%%%%%%%%%%%%%%%%%%%%%%%%
%%%%%%%%%%%%%%%%%%%%%%%%%%%%%%%%%%%%%%%%%%%%%%%%%%
\section{Markov-Dyck shifts}
%%%%%%%%%%%%%%%%%%%%%%%%%%%%%%%%%%%%%%%%%%%

Throughout this paper $N$ is a fixed positive integer larger than $1$. 
For a finite set $S$, 
we denote by 
$|S|$ its cardinality.
We consider the Dyck shift $D_N$ with alphabet 
$\Sigma = \Sigma^- \cup \Sigma^+$
where
$\Sigma^- = \{ \alpha_1,\dots,\alpha_N \},
\Sigma^+ = \{ \beta_1,\dots,\beta_N \}.
$
The symbols 
$ \alpha_i, \beta_i$
correspond to 
the brackets
$(_i,  )_i$
respectively, and have the product relations of monoid as follows:
\begin{equation}
\alpha_i \beta_j
=
\begin{cases}
 {\bold 1} & \text{ if } i=j,\\
 0  & \text{ otherwise} 
\end{cases} 
\label{eq:alphabeta}
\end{equation}
for $ i,j = 1,\dots,N$  (cf. \cite{Kr2,Kr3}).
%A word $\gamma_1\cdots\gamma_n $ of $\Sigma$
%is defined to be admissible for $D_N$ precisely if
%$\prod_{m=1}^{n} \gamma_m \ne 0.$
For a word $\omega= \omega_1 \cdots \omega_n $ of $\Sigma,$ 
we denote by $\tilde{\omega}$ its reduced form.
Namely $\tilde{\omega}$ is a word of $\Sigma \cup \{ 0, {\bold 1} \}$
obtained after applying the relations \eqref{eq:alphabeta} in $\omega$.
Then  a word $\omega$ of $\Sigma$
is said to be forbidden in $D_N$ if 
and only if $\tilde{\omega} = 0$.
Denote by 
${\frak F}_N$
the set of forbidden words.
The Dyck shift $D_N$ is defined in \cite{Kr1} 
by a subshift over $\Sigma$ whose forbidden words are ${\frak F}_N$, namely
\begin{equation}
D_N =\{ (x_n)_{n \in \Z} \in \Sigma^\Z\mid \forall k \in \Z, m\in \N, \,
(x_k,x_{k+1}, \dots,x_{k+m}) \not\in {\frak F}_N  \}. 
\end{equation}

Let $A =[A(i,j)]_{i,j=1,\dots,N}$
be an $N\times N$ matrix with entries in $\{0,1\}$.
Throughout this paper,
 $A$ is assumed to be essential 
 which means that it has no zero rows or columns.
Consider the Cuntz-Krieger algebra $\OA$ for the matrix $A$
that is the universal $C^*$-algebra generated by 
$N$ partial isometries $t_1,\dots,t_N$ subject to the following relations:
\begin{equation}
\sum_{j=1}^N t_j t_j^* = 1, 
\qquad
t_i^* t_i = \sum_{j=1}^N A(i,j) t_jt_j^* \quad \text{ for } i = 1,\dots,N
\label{eq:CK}
\end{equation}
(\cite{CK}).
Define a correspondence 
$\varphi_A :\Sigma \longrightarrow \{t_i^*, t_i \mid i=1,\dots,N\}$
by setting
$$ 
\varphi_A(\alpha_i) = t_i^*,\qquad 
\varphi_A(\beta_i) = t_i  \quad \text{ for } i=1,\dots,N.
$$
We denote by $\Sigma^*$ the set of all words 
$\gamma_1\cdots \gamma_n$ of elements of $\Sigma$.
Define the set
\begin{equation*}
{\frak F}_A = \{ \gamma_1\cdots \gamma_n \in \Sigma^* \mid
\varphi_A(\gamma_1)\cdots \varphi_A( \gamma_n) = 0 \}.
\end{equation*}
\begin{definition}
The topological Markov Dyck shift for $A$
is defined as a subshift over $\Sigma$ whose forbidden words are 
${\frak F}_A.$ 
It is written $D_A$ and called the Markov-Dyck shift for $A$ for brevity. 
\end{definition}
If $A$ is irreducible and not any permutation matrix,
the subshift $D_A$ can never be sofic (\cite[Proposition 2.1]{MathScand2011}). 
If all entries of $A$  are $1$'s, 
the $C^*$-algebra ${\mathcal{O}}_A$ 
becomes the Cuntz algebra ${\mathcal{O}}_N$
of order $N$ and
the subshift $D_A$ 
becomes the Dyck shift $D_N$
with $2N$ brackets (\cite{Cuntz}).
We note that 
 $\alpha_i \beta_j\in {\frak F}_A$  if $i\ne j$,
 and
 $\alpha_{i_n}\cdots \alpha_{i_1} \in {\frak F}_A$  
if and only if 
$\beta_{i_1}\cdots \beta_{i_n} \in {\frak F}_A$.

Let $G =(V, E)$ 
be a finite directed graph with vertex set $V$ and edge set 
$E$.
We denote by $s(e)$ the initial vertex of $e \in E$
 and by $t(e)$ the final vertex, respectively.
We assume that the cardinalities of $V$ and of $E$ are both finite 
and write
$V = \{v_1,\dots,v_{N_0}\}$ 
and 
$E = \{e_1,\dots,e_{N_1}\}$.
We also assume that each vertex of $G$
has at least one in-coming edge 
and at least one out-going edge.
The edge matrix 
$A^G =[A^G(i,j)]_{i,j=1}^{N_1}$
for $G$ is an $N_1 \times N_1$ transition matrix 
with entries in $\{0,1\}$ which is defined  by
\begin{equation}
A^G(i,j) = 
\begin{cases}
1 & \text{ if } t(e_i) = s(e_j), \\
0 & \text{ otherwise.}
\end{cases}
\end{equation}
In \cite{KMMunster}, 
we have defined the Markov-Dyck shift $D_G$ for the graph $G$ 
as the Markov-Dyck shift $D_{A^G}$ for the matrix $A^G$,
and presented formulae of the zeta function $\zeta_{D_G}(z)$
and the topological entropy $h({D_G})$.
A finite matrix $M$ with entries in $\{0,1\}$
does not necessarily arise from a finite graph as $M = A^G$.
The lemma below is  easy to prove.
For the sake of completeness, we provide its proof.
\begin{lemma}\label{lem:MAG}
Let $M =[M(i,j)]_{i,j=1}^N$
be an essential $N\times N$ matrix with entries in $\{0,1\}$.
Let us denote by
$M_r[i] =[M(i,j)]_{j=1}^N$
and
$M_c[j] =[M(i,j)]_{i=1}^N$
the $i$th row vector
and the $j$th column vector
 for $i,j=1,\dots,N$
respectively.
Then the following three conditions are equivalent:
\begin{enumerate}
\renewcommand{\theenumi}{\roman{enumi}}
\renewcommand{\labelenumi}{\textup{(\theenumi)}}
\item
There exists a finite directed graph $G$ such that
$M = A^G$.
\item
For any $i_1, i_2 \in \{1,2,\dots,N\}$,
\begin{equation}
M_r[i_1] = M_r[i_2] 
\quad
\text{ or }
\quad
\langle M_r[i_1] \mid M_r[i_2] \rangle =0. \label{eq:Mi1Mi2}
\end{equation}
\item
For any $j_1, j_2 \in \{1,2,\dots,N\}$,
\begin{equation}
M_c[j_1] = M_c[j_2] 
\quad
\text{ or }
\quad
\langle M_c[j_1] \mid M_c[j_2] \rangle =0, \label{eq:Mj1Mj2}
\end{equation}
\end{enumerate}
where 
$\langle \/ \cdot \/  \mid \/ \cdot \/ \rangle$ means the inner product
of vectors. 
\end{lemma}
\begin{proof}
(i) $\Longrightarrow$ (ii):
 Suppose that
there exists a finite directed graph $G$ such that
$M = A^G$.
For two edges $e_{i_1}, e_{i_2} \in E$,
if $t(e_{i_1}) = t(e_{i_2})$, 
then  
$M_r[i_1] = M_r[i_2]$,
otherwise
$
\langle M_r[i_1] \mid M_r[i_2] \rangle =0.
$

(iii) $\Longrightarrow$ (i):
Assume that the $N \times N$ matrix $M$
satisfies the condition
\eqref{eq:Mj1Mj2}.
We will construct a finite directed graph $G=(V,E)$ 
such that $M = A^G$ as follows.  
Define an equivalence relation 
$j_1 \sim j_2$ in $\{1,2,\dots,N\}$
by $M_c[j_1] = M_c[j_2]$.
Denote by $[j]_c$ 
the equivalence class of 
$j \in \{1,2,\dots,N\}$.
Then the vertex set $V$ is defined by 
the set of equivalence classes 
$\{ [j]_c \mid  j \in \{1,2,\dots,N\} \}$.
Define an edge labeled 
$e_i$ from 
$[i]_c$ to $[j]_c$ if $M(i,j) =1$.
If there exist edges 
from
$[i]_c$ to $[j_1]_c$ labeled $e_i$ 
and
$[i]_c$ to $[j_2]_c$ labeled $e_i$,
then 
$M(i,j_1) = M(i,j_2) =1$.
By the condition 
\eqref{eq:Mj1Mj2},
one has 
$[j_1]_c = [j_2]_c$.
Hence the labeled graph is well-defined.
Then 
as $s(e_j) = [j]_c$, 
the condition
$t(e_i) = s(e_j)$
is equivalent to
the condition
$M(i,j) =1$.
Hence we have 
$A^G = M$.

(ii) $\Longrightarrow$ (iii):
Suppose that 
there exist distinct 
$j_1 \ne j_2 \in \{1,2,\dots,N\}$
such that 
$M_c[j_1] \ne M_c[j_2]$
and
$
\langle M_c[j_1] \mid M_c[j_2] \rangle \ne 0.
$
The condition 
$M_c[j_1] \ne M_c[j_2]$
implies that 
there exists $i_1$ such that
$M(i_1,j_1) \ne M(i_1,j_2)$.
The condition
$
\langle M_c[j_1] \mid M_c[j_2] \rangle \ne 0
$
implies that
there exists $i_2$ such that
$M(i_2,j_1) = M(i_2,j_2) =1$
so that
$
\langle M_r[i_1] \mid M_r[j_2] \rangle \ne 0,
$
a contradiction to the condition (ii).
\end{proof}
The matrix 
$
\begin{bmatrix}
1 & 1 \\
1 & 0 
\end{bmatrix}
$ 
is called the Fibonacci matrix.
It can not arise from a finite directed graph
as an edge matrix.  

For a finite directed $G = (V,E)$,
we have another transition matrix $A_G$,
which is an $N_0 \times N_0$  matrix 
$A_G =[A_G(i,j)]_{i,j=1}^{N_0}$ 
defined by
\begin{equation}
A_G(i,j) = 
\begin{cases}
1 & \text{ if there exists an edge from } v_i \text{ to } v_j, \\
0 & \text{ otherwise.}
\end{cases}
\end{equation}
The matrix $A_G$ is called the vertex matrix for the graph $G$.
It has its entries in $\{0,1\}$.
\begin{definition}
Let   $G=(V,E)$ be an essential  finite directed graph.
 \begin{enumerate}
\renewcommand{\theenumi}{\roman{enumi}}
\renewcommand{\labelenumi}{\textup{(\theenumi)}}
\item 
The Markov-Dyck shift $D_{A^G}$ 
for the edge matrix $A^G$
is called the {\it Markov-Dyck shift of edge  type\/}  for $G$,
and written $D_G^E$.
\item 
The Markov-Dyck shift $D_{A_G}$ 
for the vertex matrix $A_G$
is called the {\it Markov-Dyck shift of vertex type\/} for $G$,
and written $D_G^V$.
\end{enumerate}
\end{definition}
It is obvious that any  
finite matrix $M$ with entries in $\{0,1\}$
can arise from a finite graph $G$ 
such that $M = A_G$.
By Lemma \ref{lem:MAG},
one sees that
the class of  Markov-Dyck shifts of edge type
is a subclass of  Markov-Dyck shifts of vertex type.
As is well-known that for a finite directed graph $G$
the topological Markov shift $X_{A^G}$ defined by the edge matrix $A^G$
is topologically conjugate to 
the topological Markov shift $X_{A_G}$ defined by the vertex matrix $A_G$.
The Markov-Dyck shifts however do not have this property.
Let $G_1$ be the following graph (Figure 1).
%%%%%%%%%%%%%%%%%%%%%%%%%%%%%%%%%%%%%%%%%%%%%%%%%%%%%%%%%%%%%%%%%%%%%%%%%
\begin{figure}[htbp]
\begin{center}
%WinTpicVersion3.08
\unitlength 0.1in
\begin{picture}( 25.0000,  5.8200)(  7.4000,-16.7100)
% ELLIPSE 2 0 3 0
% 4 1251 1368 1339 1466 1162 1368 1162 1368
% 
\special{pn 8}%
\special{ar 1252 1368 88 98  0.0000000 6.2831853}%
% ELLIPSE 2 0 3 0
% 4 2722 1368 2810 1466 2810 1368 2810 1368
% 
\special{pn 8}%
\special{ar 2722 1368 88 98  0.0000000 6.2831853}%
% ELLIPSE 2 0 3 0
% 4 2001 1368 1361 1065 1361 1441 2604 1441
% 
\special{pn 8}%
\special{ar 2002 1368 640 304  0.2500449 2.9054568}%
% SARROW 2 0 3 1
% 2 2617 1450 2621 1444
% 
\special{pn 8}%
\special{pa 2618 1450}%
\special{pa 2622 1444}%
\special{fp}%
\special{sh 1}%
\special{pa 2622 1444}%
\special{pa 2568 1488}%
\special{pa 2592 1488}%
\special{pa 2602 1512}%
\special{pa 2622 1444}%
\special{fp}%
% ELLIPSE 2 0 3 0
% 4 1986 1392 1346 1695 2590 1319 1346 1319
% 
\special{pn 8}%
\special{ar 1986 1392 640 304  3.3777285 6.0335374}%
% SARROW 2 0 3 1
% 2 1368 1314 1365 1320
% 
\special{pn 8}%
\special{pa 1368 1314}%
\special{pa 1366 1320}%
\special{fp}%
\special{sh 1}%
\special{pa 1366 1320}%
\special{pa 1414 1270}%
\special{pa 1390 1272}%
\special{pa 1378 1252}%
\special{pa 1366 1320}%
\special{fp}%
% ELLIPSE 2 0 3 0
% 4 949 1368 1158 1588 1136 1345 1169 1410
% 
\special{pn 8}%
\special{ar 950 1368 210 220  0.1798535 6.1660766}%
% SARROW 2 0 3 1
% 2 1152 1419 1154 1408
% 
\special{pn 8}%
\special{pa 1152 1420}%
\special{pa 1154 1408}%
\special{fp}%
\special{sh 1}%
\special{pa 1154 1408}%
\special{pa 1122 1470}%
\special{pa 1144 1460}%
\special{pa 1162 1478}%
\special{pa 1154 1408}%
\special{fp}%
% ELLIPSE 2 0 3 0
% 4 3031 1368 2822 1588 2811 1410 2844 1345
% 
\special{pn 8}%
\special{ar 3032 1368 210 220  3.2587014 6.2831853}%
\special{ar 3032 1368 210 220  0.0000000 2.9617392}%
% SARROW 2 0 3 1
% 2 2825 1331 2824 1342
% 
\special{pn 8}%
\special{pa 2826 1332}%
\special{pa 2824 1342}%
\special{fp}%
\special{sh 1}%
\special{pa 2824 1342}%
\special{pa 2850 1278}%
\special{pa 2830 1290}%
\special{pa 2810 1274}%
\special{pa 2824 1342}%
\special{fp}%
% STR 2 0 3 0
% 3 1251 1286 1251 1368 5 0
% 1
\put(12.5100,-13.6800){\makebox(0,0){1}}%
% STR 2 0 3 0
% 3 2722 1286 2722 1368 5 0
% 2
\put(27.2200,-13.6800){\makebox(0,0){2}}%
\end{picture}%
\end{center}
\caption{}
\end{figure}
%%%%%%%%%%%%%%%%%%%%%%%%%%%%%%%%%%%%%%%%%%%%%%%%%%%%%%%%%%%%%%%%%%%%%%%%%%%
The vertex matrix $A_{G_1}$ 
and the edge matrix $A^{G_1}$ are written as
\begin{equation}
A_{G_1} =
\begin{bmatrix}
1 & 1 \\
1 & 1 
\end{bmatrix},
\qquad
A^{G_1} =
\begin{bmatrix}
1 & 1 & 0 & 0 \\
1 & 1 & 0 & 0 \\
0 & 0 & 1 & 1 \\
0 & 0 & 1 & 1   
\end{bmatrix}
\end{equation}
respectively. 
Then the Markov-Dyck shift $D_{G_1}^V$ of vertex type 
is nothing but 
the Dyck shift $D_2$,
whereas 
the Markov-Dyck shift $D_{G_1}^E$ of edge type 
is not $D_2$.
 Both $D_{G_1}^V$ and $D_{G_1}^E$ have 4 fixed points as subshifts.
 The former $D_{G_1}^V$ has 4 periodic points with least period 2.
The latter $D_{G_1}^E$ has 6 periodic points with least period 2.
Hence $D_{G_1}^V$ is not topologically conjugate to $D_{G_1}^E$.

%%%%%%%%%%%%%%%%%

A Dyck $n$-path is a continuous broken directed line
on the upper half plane consisting of vectors $(1,1)$ called rise 
and $(1,-1)$ called fall.
It starts at the origin with rise  and ends at $(2n,0)$ with fall
(see \cite{Deu,GJ}, etc.). 
Let
$\gamma = (\gamma_1,\dots,\gamma_{2n})$
be a Dyck $n$-path.
Hence each $\gamma_i$ is a rise or a fall.
If $\gamma_i$ is a rise, 
there exists the smallest 
$k=1,2,\dots,2n-i$
satisfying the following two conditions:
\begin{enumerate}
\renewcommand{\theenumi}{\roman{enumi}}
\renewcommand{\labelenumi}{\textup{(\theenumi)}}
\item $\gamma_{i+k}$ is a fall.
\item $(\gamma_{i+1}, \gamma_{i+2},\dots,\gamma_{i+k-1})$ 
is a Dyck $\frac{k-1}{2}$-path
(hence $k-1$ is even), which starts at the terminal vertex of $\gamma_i$
and ends at the source vertex of $\gamma_{i+k}$.
\end{enumerate}
 We call the edge $\gamma_{i+k}$ the partner of $\gamma_i$.

Let $G=(V,E)$
be a finite directed graph. 
Denote by  $G^* =(V^*,E^*)$ 
the transposed graph of $G$.
The vertex set $V^*$ is $V$ 
and the edge set $E^*$
consists of the edges reversing its direction of the edges of $G$.
For an edge $e \in E$, we denote by $e^*$ the edge of $G^*$ 
obtained by reversing the direction of $e$, so that
$t(e^*) = s(e), s(e^*) = t(e)$ for $e \in E$. 
Recall that 
the edge set $E$ of $G$
is denoted by 
$\{e_1,\dots,e_{N_1}\}$ and 
the edge set $E^*$ of $G^*$
is written as 
$\{e_1^*,\dots,e_{N_1}^*\}$.
Put 
$ \Sigma_E^- = E^*, \Sigma_E^+ = E$
and
$
\Sigma_G^E = \Sigma_E^-\cup\Sigma_E^+.
$ 
A $G$-{\it Dyck\/} $n$-{\it path of edge type\/}
for $n=1,2,\dots$
is a Dyck $n$-path $(x_1,\dots,x_{2n})$ 
labeled 
elements of $\Sigma_G^E $
satisfying the following rules:
%\begin{enumerate}
%\item 

(1E) a rise is labeled  $e_i^*$ for some $i=1,\dots,N_1$,
%\item 

(2E)
a fall is labeled   $e_i$ for some $i=1,\dots,N_1$,
%\item 

(3E)
the partner of a rise labeled $e_i^*$ is labeled $e_i$,
%\item 

(4E)
a rise labeled  $e_i^*$ follows a rise labeled  $e_j^*$
if and only if $t(e_j^*) = s(e_i^*)$,
%\item 

(5E)
a rise labeled  $e_i^*$ follows a fall labeled  $e_j$
if and only if $t(e_j) = s(e_i^*)$,
%\item 

(6E)
a fall labeled  $e_i$ follows a fall labeled  $e_j$
if and only if $t(e_j) = s(e_i)$,
%\item 

(7E)
a fall labeled  $e_i$ follows a rise labeled  $e_j^*$
if and only if $e_j = e_i$.
%\end{enumerate}

Similarly, 
for a vertex $v \in V$, 
we denote by $v^*$ the corresponding vertex of $G^*$ 
obtained by the transposed graph $G^*=(V^*,E^*)$.
The vertex matrix $A_{G*}$ for $G^*$ 
satisfy the relations
$$
A_{G^*}(i,j) = A_G(j,i) \qquad \text{ for } i,j \in \{1,2,\dots,N_0\}.
$$
Recall that 
the vertex  set $V$ of $G$
is denoted by 
$\{v_1,\dots,v_{N_0}\}$ and 
the vertex set $V^*$ of $G^*$
is written as 
$\{v_1^*,\dots,v_{N_0}^*\}$.
Put 
$ \Sigma_V^- = V^*, \Sigma_V^+ = V$
and
$
\Sigma_G^V = \Sigma_V^-\cup\Sigma_V^+.
$ 
A $G$-{\it Dyck\/} $n$-{\it path of vertex type\/}
for $n=1,2,\dots$
is a Dyck $n$-path $(x_1,\dots,x_{2n})$ 
labeled 
elements of $\Sigma_G^V $
satisfying the following rules:
%\begin{enumerate}
%\item 

(1V)
a rise is labeled  $v_i^*$ for some $i=1,\dots,N_0$,
%\item 

(2V)
a fall is labeled   $v_i$ for some $i=1,\dots,N_0$,
%\item
 
(3V)
the partner of a rise labeled $v_i^*$ is labeled $v_i$,
%\item 

(4V)
a rise labeled  $v_i^*$ follows a rise labeled  $v_j^*$
if and only if $A_{G^*}(j,i) =1$,
%\item 

(5V)
a rise labeled  $v_i^*$ follows a fall labeled  $v_j$
if and only if $A_{G}(j,k) =A_{G^*}(k,i) =1$ for some $v_k$,
%\item 

(6V)
a fall labeled  $v_i$ follows a fall labeled  $v_j$
if and only if $A_{G}(j,i) =1$,
%\item 

(7V)
a fall labeled  $v_i$ follows a rise labeled  $v_j^*$
if and only if $v_j = v_i$.
%\end{enumerate}
%%%%%%%%%%%%%%%%%%%%%%%%%%%

The Dyck shift $D_G^E$ of edge type 
is regarded to have its symbols in $E^* \cup E$
under the identification
$\Sigma^- = E^*, \Sigma^+ = E$,
and  
the Dyck shift $D_G^V$ of vertex type 
is regarded to have its symbols in $V^* \cup V$
under the identification
$\Sigma^- = V^*, \Sigma^+ = V$.

We note the following lemma
\begin{lemma} Keep the above notations. 
\begin{enumerate}
\renewcommand{\theenumi}{\roman{enumi}}
\renewcommand{\labelenumi}{\textup{(\theenumi)}}
\item 
Any admissible word of the Dyck shift $D_G^E$ of edge type is regarded as 
a part of a labeled broken directed line of $G$-Dyck path of edge type.
Conversely a labeled broken directed line of  $G$-Dyck path of edge type
is an admissible 
word of the Dyck shift $D_G^E$ of edge type.
\item 
Any admissible word of the Dyck shift $D_G^V$ of vertex type is regarded as 
a part of a labeled broken directed line of $G$-Dyck path of vertex type.
Conversely a labeled broken directed line of $G$-Dyck path of vertex type
is an admissible 
word of the Dyck shift $D_G^V$ of vertex type.
\end{enumerate}
\end{lemma}
\begin{proof}
(i) is clear from the definition of admissible words of the Dyck shift $D_G^E$
of edge type.

(ii) 
Let $t_1,\dots,t_{N_0}$
be the partial isometries satisfying the relations \eqref{eq:CK}
for the vertex matrix $A_G$ of $G$.
For $i, j= 1,2,\dots, N_0$, we have
$\beta_j \alpha_i$ is admissible in $D_G^V$ 
if and inly if
$t_j t_i^* \ne 0$ by definition.
Since
$t_j t_i^* = t_j t_j^* t_j t_i^* t_i t_i^*$,
the condition $t_j t_i^* \ne 0$
 is equivalent to the condition that
$ t_j^* t_j t_i^* t_i \ne 0$.
As
\begin{equation*}
t_i^* t_i = \sum_{k=1}^{N_0} A_G(i,k) t_k t_k^* \qquad
t_j^* t_j = \sum_{k=1}^{N_0} A_G(j,k) t_k t_k^*,
\end{equation*}
the condition that
$ t_j^* t_j t_i^* t_i \ne 0$
is equivalent to the condition
$
 A_G(i,k) = A_G(j,k) =1 
$
for some 
$k=1,\dots,N_0$.
This shows that the condition  
$\beta_j \alpha_i$ is admissibe in $D_G^V$ 
is equivalent to the condition 
(5V)
of $G$-Dyck $n$-path of vertex type.
It is direct to see that the other conditions
(1V), (2V), (3V), (4V), (6V), (7V)
are compatible to the definitions of giving admissible words 
of the Dyck shift $D_G^V$ of vertex type. 
 \end{proof}

We remark that 
a finite path of vertices of 
a labeled broken directed line  of the $G$-Dyck path of edge type is 
not necessarily an admissible 
word of the Dyck shift $D_G^V$ of vertex type.
Consider the following correspondences in $G$-Dyck paths:
\begin{equation}
\begin{cases}
\text{ a fall } e \in E  
& \longrightarrow 
\text{ the source } s(e) \in V \text{ of } e, \\ 
\text{ a rise } e^* \in E^*  
& \longrightarrow 
\text{ the terminal  } t(e^*) \in V^* \text{ of } e^*. 
\end{cases} \label{eq:EV}
\end{equation}
The rules $(1E),\dots,(7E)$ and
$(1V),\dots,(7V)$ ensure us the following lemma.
\begin{lemma} Keep the above notations. 
\begin{enumerate}
\renewcommand{\theenumi}{\roman{enumi}}
\renewcommand{\labelenumi}{\textup{(\theenumi)}}
\item 
Any sequence of vertices of a $G$-Dyck $n$-path of edge type 
yields a labeled sequence by $\Sigma_G^V$ 
of a $G$-Dyck $n$-path of vertex type 
by the correspondence \eqref{eq:EV}.
\item 
Any labeled sequence by $\Sigma_G^V$ of a $G$-Dyck $n$-path of vertex type
is realized as a sequence of vertices of a $G$-Dyck $n$-path of edge type
by the correspondence \eqref{eq:EV}.
\end{enumerate}
\end{lemma}
By the above lemma, it is reasonable to  
define a $1$-block map 
$\varPhi:E \cup E^* \longrightarrow V \cup V^*$ 
by
\begin{equation*}
\begin{cases}
\varPhi(e) & = s(e) \in V \quad \text{ for } e \in E, \\
\varPhi(e^*) & = t(e^*)(=s(e)) \in V^* \quad \text{ for } e^* \in E^*
\end{cases}
\end{equation*}
Hence we have 
\begin{proposition}\label{prop:oneblock}
The $1$-block map $\varPhi:E \cup E^* \longrightarrow V \cup V^*$ 
induces a factor code $\varphi =\varPhi_\infty: D_G^E \longrightarrow D_G^V$.
\end{proposition}
For $e_{i,k} \in E$ with 
$s(e_{i,k}) = v_i \in V$ and $t(e_{i,k}) = v_k \in V$,
and
$e_{k,j}^* \in E^*$ with 
$s(e_{k,j}^*) = v_k^* \in V^*$ 
and $t(e_{k,j}^*) = v_j^* \in V^*$,
then 
the word $(e_{i,k}, e_{k,j}^*)$ is admissible in $D_G^E$
and
the word $(v_{i}, v_{j}^*)$ is admissible in $D_G^V$
such that
\begin{equation*}
\varPhi\left(\overset{v_i}{\qquad\underset{\qquad{v_k}}{\searrow}}
\overset{\quad{v_j^*}}{\nearrow} \qquad \right) 
=(i\searrow \nearrow j^*), 
\qquad
\varPhi(e_{i,k}, e_{k,j}^*) =(v_{i}, v_{j}^*).
\end{equation*}
%\begin{equation*}
%\varPhi\overset{\qquad{v_i}}{(}\!\!\underset{\qquad{v_k}}{\searrow}
%\overset{\quad{v_j^*}}{\nearrow} ) =(i\searrow \nearrow j^*), 
%\qquad \varPhi(e_{i,k}, e_{k,j}^*) =(v_{i}, v_{j}^*).
%\end{equation*}
In the above situation, 
we call the vertex 
$v_k(=v_k^*)$ a valley. 
Hence the factor map
$\varphi:D_G^E \longrightarrow D_G^V$ erases the valleys. 
We will show that the factor map $\varphi$ is finite-to-one,
so that the equality of the topological entropy 
$h_{top}(D_G^E) = h_{top}(D_G^V)$ holds.

We provide  the height functions on $D_G^E$.
These functions on the Dyck shift $D_N$
have been first introduced by W. Krieger in \cite{Kr1}.
For 
$x = (x_n)_{n \in {\mathbb{Z}}} \in D_G^E$,
we set the height function 
\begin{align*}
H_0(x) & = 0, \\
H_m(x) & = \sum_{k=0}^{m-1}
(\chi_-(x_k)- \chi_+(x_k)), \qquad m \in {\mathbb{N}}, \\
H_{-m}(x) & = \sum_{k=-1}^{-m}
(- \chi_-(x_k) + \chi_+(x_k)), \qquad m \in {\mathbb{N}}
\end{align*} 
where
\begin{equation*}
\chi_-(x_k)
=
{\begin{cases}
1 & \text{ if } x_k \in \Sigma^-, \\
0 & \text{ if } x_k \in \Sigma^+, 
\end{cases}
} \qquad
\chi_+(x_k)
=
{\begin{cases}
0 & \text{ if } x_k \in \Sigma^-, \\
1 & \text{ if } x_k \in \Sigma^+. 
\end{cases}
}
\end{equation*}
\begin{definition}
For $x =(x_n)_{n \in \Z} \in D_G^E$,
\begin{enumerate}
\renewcommand{\theenumi}{\roman{enumi}}
\renewcommand{\labelenumi}{\textup{(\theenumi)}}
\item
a vertex $t(x_{m-1})(=s(x_m))$ is called a relative minimum in $x$ 
if
$x_{m-1} \in E$ and $x_m \in E^*$. 
\item
a vertex $t(x_{m-1})(=s(x_m))$ is called a minimum in $x$ if
$H_m(x) \le H_n(x) $ for all $n \in \Z$. 
\end{enumerate}
\end{definition}
\begin{lemma}
For $x =(x_n)_{n \in \Z} \in D_G^E$,
\begin{enumerate}
\renewcommand{\theenumi}{\roman{enumi}}
\renewcommand{\labelenumi}{\textup{(\theenumi)}}
\item
if a vertex $t(x_{m-1})(=s(x_m))$ is not a relative minimum in $x$,
the word  $(\varPhi(x_{m-1}),\varPhi(x_{m}))$ in $D_G^V$ 
uniquely determines the vertex $t(x_{m-1})$,
\item
if a vertex $t(x_{m-1})(=s(x_m))$ is not minimum in $x$,
the sequence $\varphi(x) \in D_G^V$ 
uniquely determines the vertex $t(x_{m-1})$,
 \item
if two vertices $t(x_{n-1})$ and $t(x_{m-1})$ are both  minimum in $x$,
then $t(x_{n-1}) =t(x_{m-1})$.
\end{enumerate}
\end{lemma}
\begin{proof}
(i)
Since  the vertex $t(x_{m-1})(=s(x_m))$ is not a relative minimum in $x$,
we have two cases.

Case 1: $x_{m-1} \in E^*$.

Since $\varPhi(x_{m-1}) $ is in $V^*$, 
we take a vertex $v_i \in V$ such that 
$\varPhi(x_{m-1})= v_i^*$. We then have
$t(x_{m-1}) = v_i^*$.

Case 2: $x_{m-1} \in E$.

The condition that the vertex $t(x_{m-1})(=s(x_m))$ is not a relative minimum in $x$
implies that
$x_m$ belongs to $E$, so that 
$\varPhi(x_{m})= v_j\in V$
for some $j$.
We then have
$t(x_{m-1}) =s(x_m) = v_j.$

(ii)
Suppose that the  vertex $t(x_{m-1})(=s(x_m))$ is not minimum in $x$.
If $t(x_{m-1})$ is not a relative minimum in $x$,
the above discussion implies that 
the word  $(\varPhi(x_{m-1}),\varPhi(x_{m}))$ in $D_G^V$ 
uniquely determines the vertex $t(x_{m-1})$.
Hence we may assume that  
$t(x_{m-1})$ is a relative minimum in $x$.
Since $t(x_{m-1})(=s(x_m))$ is not minimum in $x$,
there exists $i \in \Z$ such that
$H_i(x) < H_m(x)$.
We have two cases.

Case 1: $i >m$.

There exists $k\in \Z$ with  $m < k <i$
such that
$x_{k-1}, x_k \in E$, 
and
$H_m(x) = H_k(x)$.
We take a vertex $v_j \in V$ such that 
$\varPhi(x_{k})= v_j$. We then have
$t(x_{m-1})= t(x_{k-1}) = v_j$.

Case 2: $i < m$.

There exists $l\in \Z$ with  $i < l <m$
such that
$x_{l-1}, x_l \in E^*$, 
and
$H_m(x) = H_l(x)$.
We take a vertex $v_j \in V$ such that 
$\varPhi(x_{l-1})= v_j$. We then have
$t(x_{m-1})= t(x_{l-1}) = v_j$.

(iii)
Suppose that two vertices $t(x_{n-1})$ and $t(x_{m-1})$ are both  minimum in $x$,
so that $H_n(x) = H_m(x)$.
Assume that $n <m$.
The word $(x_n,x_{n+1}, \dots, x_{m-1})$ is a $G$-Dyck path of edge type
so that the vertices 
$s(x_n)$ and $t(x_{m-1})$ are the same.
This implies that $t(x_{n-1}) =t(x_{m-1})$.
\end{proof}

\begin{proposition}\label{prop:factor}
Suppose that $G$ does not have multi-edges.
Let $\varphi:D_G^E \longrightarrow D_G^V$ be the factor code defined
in Proposition \ref{prop:oneblock}.
For $x =(x_n)_{n \in \Z} \in D_G^E$, we have
\begin{enumerate}
\renewcommand{\theenumi}{\roman{enumi}}
\renewcommand{\labelenumi}{\textup{(\theenumi)}}
\item
if $x$ does not have a minimum vertex,
then $\varphi$ is injective at $x$, that is,
\begin{equation*}
\varphi^{-1}(\varphi(x)) = x,
\end{equation*}
\item
if $x$  has a minimum vertex,
then 
\begin{equation*}
| \varphi^{-1}(\varphi(x))| \le N_0 = |V|.
\end{equation*}
\end{enumerate}
Therefore $\varphi:D_G^E \longrightarrow D_G^V$ is a finite-to-one  
factor code.
\end{proposition}
\begin{proof}
(i)
Suppose that $x = (x_n)_{n \in \Z}$ does not have a miniumum vertex.
By (ii) of the above lemma, the sequence
$\varphi(x)$ determines the sequence $t(x_n), n \in \Z$ of vertices.
Each symbol $x_n$ is an edge of $E$ or of $E^*$, and an edge is determined by
the vertices $t(x_n), t(x_{n-1})(=s(x_n))$,
 so that the code $\varphi$ is injective at $x$.

(ii)
Suppose that $x$ has a minimum vertex at $t(x_{m-1})$ for
some $m \in \Z$.
Then  
the vertex  $t(x_{m-1})$ is a valley and
$x_{m-1} \in E, x_m \in E^*$.
By (iii) of the above lemma,
other minimum vertices are the same as 
the vertex 
$t(x_{m-1})$.
Hence we have 
\begin{align*}
| \varphi^{-1}(\varphi(x))| 
& = | \{ k \in \{1,2,\dots,N_0\} \mid 
A_G(s(x_{m-1}), k)= A_{G^*}(k,t(x_{m})) =1 \} |\\
&  \le  N_0 = |V|.
\end{align*}
\end{proof} 
\begin{theorem}\label{thm:entropy}
Suppose that $G$ does not have multi-edges.
We then have
$h_{top}(D_G^V) = h_{top}(D_G^E)$.
\end{theorem}
\begin{proof}
Since there exists a factor code
 $\varphi:D_G^E \longrightarrow D_G^V$,
 the inequality
$h_{top}(D_G^V) \le h_{top}(D_G^E)$
 is clear.
The $1$-block map  $\varPhi$ naturally induces a map
 $\varPhi_*:B_*(D_G^E) \longrightarrow B_*(D_G^V)$
 between admissible words.
 It is not necessarily one-to-one at minimal points of words.
We then have
\begin{equation*}
| B_n(D_G^E)|\le N_0 \cdot | B_n(D_G^V)|, \qquad n \in \N
\end{equation*}  
Therefore we have
$h_{top}(D_G^E) \le h_{top}(D_G^V)$.
\end{proof}
Concerning embedding of the Markov-Dyck shifts,
we have the following proposition.
\begin{proposition}
Suppose that $G$ does not have multi-edges.
There exists an embedding of $D_G^E$ 
into the 3rd power shift of $D_G^V$.
\end{proposition}
\begin{proof}
Let
$t_i,i=1,\dots,N_0$
be partial isometries satisfying the relations  
\eqref{eq:CK} for the vertex matrix $A_G$.
For an edge $e_n \in E$ with $s(e_n) = v_i, t(e_n) = v_j$,
define a partial isometry
$S_n = t_i t_jt_j^*$.
It is easy to see
that the family $S_1,\dots, S_{N_1}$ satisfies the relations
\eqref{eq:CK} for the edge matrix $A^G$,
This implies that the
 correspondence 
$\Psi: E\cup E^* \longrightarrow (V \cup V^*)^{[3]}$
defined by
\begin{equation*}
\Psi(e_n) = (v_i, v_j, v_j^*),\qquad
\Psi(e_n^*) = (v_j, v_j^*, v_i^*)
\end{equation*}  
induces an embedding
of $D_G^E$ into the 3rd power shift $(D_G^V)^{[3]}$
of $D_G^V$.
\end{proof}

%%%%%%%%%%%%%%%%%%%%%%%%%%%%
\section{The zeta functions of Mrkov-Dyck shifts of vertex type}
%%%%%%%%%%%%%%%%%%%%%%%%%%%%%%%%%%%%%%%%%%%%%%%%%%%%%%%%%%%%%
In what follows, we fix an arbitrary $N \times N$ matrix 
$A =[A(i,j)]_{i,j=1}^N $ with entries in $\{0,1\}$.
We will study the Markov-Dyck shift 
$D_A$ and present a formula of the zeta function
$\zeta_{D_A}(z)$.
In \cite{KMMunster}, 
a formula of the zeta function of the Markov-Dyck shifts
of edge type has been presented.
The Markov-Dyck shifts of edge type form a subclass of the class of 
Markov-Dyck shifts. 
In this section, 
we will study general  Markov-Dyck shift 
$D_A$ and present a formula of its zeta function
$\zeta_{D_A}(z)$.    
For the $N \times N$ matrix $A$, 
let
$v_1,\dots,v_N$ be $N$-vertices.
Define a directed edge from $v_i$ to $v_j$ if
$A(i,j) = 1$.
We then have a finite directed graph written $G = (V,E)$
such that its vertex matrix $A_G$ coincides with the original matrix
$A$.

Throughout this section, 
we identify $\alpha_i$ with $v_i^*$
and
$\beta_i$ with $v_i$
for $i=1,\dots,N$, respectively.
Let
$w=(w_1,\dots,w_{2n})$
be a $G$-Dyck $n$-path of vertex type.
As in 
\cite{MatsumotoIJM},
$w$ is called a $G$-Catalan word
and satisfies the following conditions:
\begin{align*}
 \sum_{k=1}^{m} &
(\chi_-(w_k)- \chi_+(w_k)) \ge 0 \qquad \text{ for all }
m=1,2,\dots,2n\\
\intertext{ and }
 \sum_{k=1}^{2n} &
(\chi_-(w_k) - \chi_+(w_k)) =0.
\end{align*} 
Denote by
 $C_n^A$
 the set of $G$-Dyck $n$-pathes of vertex type.
For $i=1,\dots,N$, 
put
\begin{equation*}
C_n^A(i) = \{ (w_1,\dots,w_{2n}) \in C_n^A \mid
(\alpha_i,w_1,\dots,w_{2n},\beta_i) \in C_{n+1}^A
\}.
\end{equation*}
Denote by $c_n^A(i)$ the cardinarity
$|C_n^A(i)| $ of the set $C_n^A(i). $
We set $c_0^A(i) =1$.
Combinatorial properties of the sequence
 $c_n^A(i), n=0,1,\dots$
 have been studied in
 \cite[Section 4]{MatsumotoIJM}. 
For $i=1,\dots,N$, let
$f_i^A(z)$ be the generating function of the sequence 
$c_n^A(i), n=0, 1,2,\dots:$
\begin{equation*}
f_i^A(z) = \sum_{n=0}^\infty c_n^A(i)z^n.
\end{equation*}
Since one knows (\cite[Section 4]{MatsumotoIJM})
\begin{equation*}
C_{n+1}^A(i)= \bigcup_{k=0}^n 
              \bigcup_{\substack{j\\
                         A(j,i)=1}} 
               C_k^A(j) \times C_{n-k}^A(i), 
\end{equation*}
we have
\begin{equation*}
c_{n+1}^A(i)= \sum_{k=0}^n \sum_{j=1}^N A(j,i) c_k^A(j) c_{n-k}^A(i), 
\end{equation*}
so that the identity
\begin{equation} 
f_i^A(z)=  1 + z f_i^A(z) \sum_{j=1}^N A(j,i) f_j^A(z) \label{eq:fia}
\end{equation}
holds (\cite[Proposition 4.2]{MatsumotoIJM}).
Let $X_A$ be the shift space over $\Sigma^+ =V$
of the topological Markov shift
defined by the matrix $A$:
\begin{equation*}
X_A =
\{ (x_n)_{n \in \Z}\in (\Sigma^+)^{\Z} \mid A(x_n,x_{n+1}) =1 
\text{ for all } n \in \Z\}.
\end{equation*}
For $n, k \in \N$, we set
\begin{align*}
C_{n,k}^{A,+} = \{
&(w_1,\dots,w_{2n},\beta_{i_1},\dots,\beta_{i_k}) \in B_{2n+k}(D_A)\mid \\
&(w_1,\dots,w_{2n})  \in C_n^A,\,
(\beta_{i_1},\dots,\beta_{i_k}) \in B_{k}(X_A)
\}.
\end{align*} 
For
$(w_1,\dots,w_{2n},\beta_{i_1},\dots,\beta_{i_k}) \in C_{n,k}^{A,+}$,
we set
\begin{align*}
s((w_1,\dots,w_{2n},\beta_{i_1},\dots,\beta_{i_k}))
& = \beta_{i_1}, \\
t((w_1,\dots,w_{2n},\beta_{i_1},\dots,\beta_{i_k}))
& = \beta_{i_k}.
\end{align*} 
We put
$$
C_A^+ = \bigcup_{n=1}^\infty\bigcup_{k=1}^\infty C_{n,k}^{A,+}.
$$
We then see the following lemma.
\begin{lemma}
For $\mu, \nu \in C_A^+ $,
the word $\mu \nu $ is admissible in $D_A$
if and only if
$A(t(\mu),s(\nu)) =1$.
\end{lemma}
Put $I = \{1,\dots,N\} \times \{1,\dots,N\}$.
Define an $I \times I$ matrix 
$\tilde{A}=[\tilde{A}((i,j),(k,l))]_{(i,j),(k,l)\in I}$
by
\begin{equation*}
\tilde{A}((i,j),(k,l)) = A(j,k)
\end{equation*}
and a map
$r: C_A^+ \longrightarrow I$
by
\begin{equation*}
r((w_1,\dots,w_{2n},\beta_{i_1},\dots,\beta_{i_k}))
= (\beta_{i_1},\beta_{i_k}) \in I.
\end{equation*}
Then the quadruplet 
${\mathcal{C}}_A^+
=(C_A^+, I, \tilde{A}, r)$ is a circular Markov code in the sense of Keller 
\cite{Keller}.
We then associate the following shift-invariant subset 
$\Omega_{{\mathcal{C}}_A^+}$ 
by
\begin{align*}
\Omega_{{\mathcal{C}}_A^+}
=\{& x=(x_n)_{n \in \Z} \mid
\text{ there are } \dots k_{-1} < k_0\le 0 < k_1 < \dots \text{ in } \Z \\
 &\text{ such that }
 x_{[k_i,k_{i+1})} \in C_A^+ \text{ and } 
\tilde{A}{(r(x_{[k_{i-1},k_{i})}), r(x_{[k_i,k_{i+1})})}) =1
\} \qquad (\cite{Keller}).
\end{align*}
The zeta function $\zeta(\Omega_{{\mathcal{C}}_A^+},z)$
for a shift-invariant set $\Omega_{{\mathcal{C}}_A^+}$
is similarly defined to \eqref{eq:zeta} by using a sequence of cardinalities
of periodic points of $\Omega_{{\mathcal{C}}_A^+}$. 
Following Keller \cite{Keller},
define a sequence 
$D({\mathcal{C}}_A^+,m) = \diag[d_{(i,j),(i,j)}({\mathcal{C}}_A^+,m)],
3\le m \in \N$
of 
$I \times I$-diagonal matrices
with diagonal entries 
$d_{(i,j),(i,j)}({\mathcal{C}}_A^+,m), (i,j)\in I$ 
 by
\begin{align*}
d_{(i,j),(i,j)}({\mathcal{C}}_A^+,m)
&= |\{(w_1,\dots,w_{2n},\beta_{i_1},\dots,\beta_{i_k})\in C_A^+
\mid i_i =i, i_k = j \}| \\
(&=c_n^A(i)A^{k-1}(i,j))
\end{align*}
for $m = 2n+k$,
and  a matrix-valued generating function
$F({\mathcal{C}}_A^+,z)$ by
\begin{equation*}
F({\mathcal{C}}_A^+,z) = \sum_{m=1}^\infty D({\mathcal{C}}_A^+,m) 
\tilde{A}z^m.
\end{equation*}
Denote by
$I_{N^2}$ the identity matrix of size $N^2$.
By using \cite[Theorem 1]{Keller}, we have
\begin{proposition}
$
\zeta(\Omega_{{\mathcal{C}}_A^+},z)
= \det(I_{N^2} - F({\mathcal{C}}_A^+,z))
$
\end{proposition}
We then have for $(i,j), (p,q) \in I$
\begin{align*}
F({\mathcal{C}}_A^+,z)((i,j),(p,q))
& = \sum_{m=1}^\infty D({\mathcal{C}}_A^+,m) \tilde{A}z^m ((i,j),(p,q))\\
& = \sum_{m=1}^\infty 
    \sum_{\substack{n,k\\
                    2n+k =m}}
    D({\mathcal{C}}_A^+,2n+k) \tilde{A}z^{2n+k} ((i,j),(p,q))\\
& = \sum_{n=1}^\infty 
\sum_{k=1}^\infty c_n^A(i) A^{k-1}(i,j)\tilde{A}((i,j),(p,q))z^{2n+k} \\
& = \sum_{n=1}^\infty c_n^A(i) z^{2n}
\sum_{k=1}^\infty  A^{k-1}(i,j) z^k \tilde{A}((i,j),(p,q))\\
& = (f_i^A(z^2) -1) z
\sum_{l=0}^\infty (z A)^{l}(i,j)A(j,p)\\
& = (f_i^A(z^2) -1) z (1_N - zA)^{-1}(i,j)\cdot A(j,p).
\end{align*}
We define $N \times N$ matrices 
$F^A = [F^A(i,j)]_{i,j=1}^N$
and
$H(C_A^+,z) $
by
\begin{equation*}
F^A(i,j) = (f_i^A(z^2) -1) z (1_N - zA)^{-1}(i,j)
\quad 
\text{ and }
\quad
H(C_A^+,z) = F^A\cdot A
\end{equation*}
so that
\begin{equation*}
F({\mathcal{C}}_A^+,z)((i,j),(p,q))
=
F^A(i,j) A(j,p)
\quad
\text{ and }
\quad
H(C_A^+,z)(i,p)  = 
\sum_{j=1}^N F({\mathcal{C}}_A^+,z)((i,j),(p,q)).
\end{equation*}
\begin{lemma}
$
\det(I_{N^2} -F({\mathcal{C}}_A^+,z))
=
\det(I_N -H(C_A^+,z)).
$
\end{lemma}
\begin{proof}
Let
$
U = [U((i,j),(p,q))]_{(i,j),(p,q) \in I}$
and
$V = [V((i,j),(p,q))]_{(i,j),(p,q) \in I}$
be $I \times I$ matrices 
defined by
\begin{align*}
U((i,j),(p,q))
& =
{\begin{cases}
1 & \text{ if } (i,j)=(p,q),\\
1 & \text{ if } i=p, \, j=N,\\
0 & \text{ otherwise},
\end{cases} }\\
V((i,j),(p,q))
& =
{\begin{cases}
1 & \text{ if } (i,j)=(p,q),\\
-1& \text{ if } i=p, \, j=N, \, q<N, \\
0 & \text{ otherwise}.
\end{cases}}
\end{align*}
The matrix
$(I_{N^2} -F({\mathcal{C}}_A^+,z))V$
is obtained from
$(I_{N^2} -F({\mathcal{C}}_A^+,z))$
by adding the minus of the 
$(i,N)$th column to the 
$(i,j)$th column 
for all $j=1,2,\dots,N-1$ and $i=1,2,\dots,N$,
and 
the matrix
$U(I_{N^2} -F({\mathcal{C}}_A^+,z))V$
is obtained from
$(I_{N^2} -F({\mathcal{C}}_A^+,z))V$
by adding the 
$(i,j)$th rows  
to the $(i,N)$th row
for all $j=1,2,\dots,N-1$ and $i=1,2,\dots,N$.
Hence we see 
\begin{equation*}
U(I_{N^2} -F({\mathcal{C}}_A^+,z))V((i,j),(p,q)) 
=
\begin{cases}
1 & \text{ if } (i,j)=(p,q),\, q<N, \\
0 & \text{ if } (i,j) \ne (p,q),\, q<N,\\
1- \sum_{k=1}^N F^A(i,k)A(k,p)  & \text{ if } (i,j)=(p,q),\, q =N,\\
-F^A(i,j)A(j,p) & \text{ if } j<N, \,  q=N,\\
0 & \text{ otherwise.}
\end{cases}
\end{equation*} 
Each $(p,q)$th column for $q<N$ of
the matrix $U(I_{N^2} -F({\mathcal{C}}_A^+,z))V$
has $1$ on diagonal and zero elsewhere.
Since
$$
1- \sum_{k=1}^N F^A(i,k)A(k,p) = 
1 - H(C_A^+,z)(i,p),
$$
by expanding the matrix
$U(I_{N^2} -F({\mathcal{C}}_A^+,z))V $
along the $(p,q)$th columns for 
$p=1,2,\dots,N$ with $q<N$,
we have
\begin{equation*}
\det(U(I_{N^2} -F({\mathcal{C}}_A^+,z))V )
=
\det(I_N -H(C_A^+,z)).
\end{equation*}
As
$\det(U) = \det(V) =1$,
we get the desired equality.
\end{proof}
Therefore we have
\begin{proposition}\label{prop:CAplus}
\begin{equation}
\zeta(\Omega_{{\mathcal{C}}_A^+},z)
= \frac{\det(I_N - zA)}{
\det(I_N - \diag[f_1^A(z^2),\dots,f_N^A(z^2)] zA)}.\label{eq:zetacaplus}
\end{equation}
\end{proposition}
\begin{proof}
Since
\begin{equation*}
H(C_A^+,z)
=
\diag[f_1^A(z^2)-1,\dots,f_N^A(z^2)-1] zA (I_N - z A)^{-1},
\end{equation*}
we have
\begin{align*} 
I_N - H(C_A^+,z)
& = I_N  
   - \diag[f_1^A(z^2),\dots,f_N^A(z^2)] zA (I_N - z A)^{-1}
   + zA (I_N - z A)^{-1} \\
& = (I_N - z A)^{-1} - \diag[f_1^A(z^2),\dots,f_N^A(z^2)]zA(I_N - z A)^{-1}\\
& = (I_N - \diag[f_1^A(z^2),\dots,f_N^A(z^2)]zA ) (I_N - z A)^{-1}
\end{align*}
so that the desired equality holds.
\end{proof}

\medskip

For $j \in \{1,2,\dots,N\}$ with 
$A(i,j) =1$, we put
\begin{equation*}
C_n^A[i;\{j\}] =
\{ (\alpha_i, w_1,\dots,w_{2n-2},\beta_i) \in C_n^A(j) 
\mid
(w_1,\dots,w_{2n-2}) \in C_{n-1}^A(i) \}
\end{equation*}
and
\begin{equation*}
C_n^A[j] = \bigcup^N_{\substack{i=1\\
                       A(i,j)=1}}
C_n^A[i;\{j\}],
\qquad 
C^A[j] = \bigcup_{n=1}^\infty C_n^A[j]. 
\end{equation*}
We set $c_n^A[j] = |C_n^A[j]|$.
As
$|C_n^A[i;\{j\}]| = c_{n-1}^A(i)$ if 
$A(i,j) =1$,
we have
\begin{equation}
c_n^A[j]  
=\sum_{i=1}^N A(i,j) c_{n-1}^A(i). \label{eq:cnj}
\end{equation}
Similarly for a subset
$\{ j_1,\dots,j_k \} \subset \{1,2,\dots,N\}$
with 
$A(i,j_1) =\cdots =A(i,j_k) =1$, we put
\begin{equation*}
C_n^A[i;\{ j_1,\dots,j_k \}] =
\bigcap_{m=1}^k  C_n^A[i;\{j_m\}]
\end{equation*}
and
\begin{align*}
C_n^A[\{ j_1,\dots,j_k \}] 
 &= \bigcup^N_{\substack{i=1\\
    A(i,j_1)=\cdots =A(i,j_k)=1}
    }
     C_n^A[i;\{ j_1,\dots,j_k \}],\\
C^A[\{ j_1,\dots,j_k \}] 
& = \bigcup_{n=1}^\infty C_n^A[\{ j_1,\dots,j_k \}]. 
\end{align*}
We set $c_n^A[\{ j_1,\dots,j_k \}] = |C_n^A[\{ j_1,\dots,j_k \}]|$ 
so that
\begin{equation}
c_n^A[\{ j_1,\dots,j_k \}] 
=\sum_{i=1}^N A(i,j_1)\cdots A(i,j_k) c_{n-1}^A(i). \label{eq:cnk} 
\end{equation}
%For a subset $ \{ j_1,\dots,j_k \} \subset \{1,2,\dots,N\}$,
%if there is no $i \in \{1,2,\dots,N\}$ satisfying
%$A(i,j_1)=\cdots =A(i,j_k)=1$,
%we set $C^A[\{ j_1,\dots,j_k \}] =\emptyset$ 
%and hence 
%$c_n^A[\{ j_1,\dots,j_k \}] =0$.
For a subset
$ \{ j_1,\dots,j_k \} \subset \{1,2,\dots,N\}$
if
there exists $i \in \{1,2,\dots,N\}$
such that 
$A(i,j_1)=\cdots =A(i,j_k)=1$,
we call the set
$ 
C^A[\{ j_1,\dots,j_k \}] 
$ 
the {\it Markov-Dyck code with support} 
$\{ j_1,\dots,j_k \}$.
It is easy to see that the set
$ 
C^A[\{ j_1,\dots,j_k \}] 
$ 
is a circular code.
Denote by
$ 
C^A[\{ j_1,\dots,j_k \}]^{\infty} 
$ 
the set of all two-sided sequences of alphabet
$\Sigma =\Sigma^- \cup \Sigma^+$
consisting  of free concatenations of words of 
$ 
C^A[\{ j_1,\dots,j_k \}]. 
$ 
Let $g_{C^A[\{ j_1,\dots,j_k \}]}(z) $ be the generating function
for the sequence
$c_n^A[\{ j_1,\dots,j_k \}], n=1,2,\dots$ defined by
\begin{equation*}
g_{C^A[\{ j_1,\dots,j_k \}]}(z) 
= \sum_{n=1}^\infty
c_n^A[\{ j_1,\dots,j_k \}]z^{2n}.
\end{equation*}
\begin{lemma}
\begin{enumerate}
\renewcommand{\theenumi}{\roman{enumi}}
\renewcommand{\labelenumi}{\textup{(\theenumi)}}
\item
The generating function
$g_{C^A[\{ j_1,\dots,j_k \}]}(z)$
satisfies
\begin{equation}
g_{C^A[\{ j_1,\dots,j_k \}]}(z) =
z^2 \sum_{i=1}^N A(i,j_1)A(i.j_2)\cdots A(i,j_k) f_i^A(z^2).
\end{equation}
\item
The zeta function 
$\zeta(C^A[\{ j_1,\dots,j_k \}]^{\infty}, z)$
of the shift-invariant set

$C^A[\{ j_1,\dots,j_k \}]^{\infty}\subset \Sigma^{\Z}$ 
is
\begin{equation}
\zeta(C^A[\{ j_1,\dots,j_k \}]^{\infty}, z) = 
\frac{1}{1-g_{C^A[\{ j_1,\dots,j_k \}]}(z)}.
\end{equation}
In particular for $j \in \{1,2,\dots,N\}$, we have
\begin{equation}
\zeta(C^A[\{ j \} ]^{\infty}, z) = 
\frac{1}{1-g_{C^A[\{ j \}]}(z)} = f_{j}^A(z^2).
\end{equation} 
\end{enumerate}
\end{lemma}
\begin{proof}
(i)
By \eqref{eq:cnk}, we have
\begin{align*}
g_{C^A[\{ j_1,\dots,j_k \}]}(z) 
%& = \sum_{n=1}^\infty
%    c_n^A[\{ j_1,\dots,j_k \}]z^{2n}\\
& = \sum_{n=1}^\infty
    \sum_{i=1}^N A(i,j_1)\cdots A(i,j_k) c_{n-1}^A(i) z^{2n}\\
& = z^2 \sum_{i=1}^N A(i,j_1)\cdots A(i,j_k) 
        \sum_{n=1}^\infty c_{n-1}^A(i) z^{2(n-1)}\\
& = z^2 \sum_{i=1}^N A(i,j_1)\cdots A(i,j_k) f_i^A(z^2). 
\end{align*}

(ii)
The set $C^A[\{ j_1,\dots,j_k \}]$
is a circular code, and 
the set 
$ 
C^A[\{ j_1,\dots,j_k \}]^{\infty} 
$ 
consisting of the two-sided sequences of free concatenations of words of 
$ 
C^A[\{ j_1,\dots,j_k \}]. 
$ 
Hence a well-known theorem of combinatorics 
(cf. \cite[Proposition 4.7.11]{Stanley})
ensures us the equality
\begin{equation*}
\zeta(C^A[\{ j_1,\dots,j_k \}]^{\infty}, z) 
= 
\frac{1}{1-g_{C^A[\{ j_1,\dots,j_k \}]}(z)}.
\end{equation*}
In particular
we have
\begin{equation*}
g_{C^A[\{ j \}]}(z)
= z^2 \sum_{i=1}^N A(i,j) f_i^A(z^2) 
= \frac{f_j^A(z^2) -1}{f_j^A(z^2) }
= 1 - \frac{1}{f_j^A(z^2)}
\end{equation*}
so that
\begin{equation}
\zeta(C^A[\{ j \} ]^{\infty}, z)  
=\frac{1}{1 - g_{C^A[\{ j \}]}(z)} 
=f_j^A(z^2).
\end{equation}
\end{proof}
We call a subset
$ \{ j_1,\dots,j_k \} \subset \{1,2,\dots,N\}$
 {\it a support subset\/} if for any $i \in \{1,2,\dots,N\}$
 there exists $l=1,\dots,k$ such that 
$A(i,j_l)=1$.
The set $\{1,2,\dots,N\}$ itself is a support subset. 
For a shift-invariant subset $C$ of $D_A$,
denote by 
$P_n(C)$ 
the set of $n$-periodic points of $C$.
We set
\begin{equation}
{C^A}^{\infty} =\bigcup_{ \{ j_1,\dots,j_k \} \subset \{1,\dots,N\}} 
C^A[\{ j_1,\dots,j_k \}]^{\infty} \quad \subset \Sigma^{\Z}.
\end{equation}
By the principle of inclusion of exclusion in combinatorics
(cf. \cite[2.1]{Stanley}),
we have
\begin{lemma}
Let $J =\{ j_1,\dots,j_k \}$ be a support subset of 
$\{1,2,\dots,N\}$. Then we have
\begin{align*}
&P_n({C^A}^{\infty})\\ 
= & \bigcup_{l=1}^k P_n(C^A[\{ j_l \} ]^{\infty})
-
 \bigcup_{\{ j_1,j_2\}\subset J}P_n(C^A[\{ j_1,j_2 \}]^{\infty}) \\
 \cdots
& {{(-1)}^{m+1}} \bigcup_{\{j_1,\dots,j_m\} \subset J}
P_n(C^A[\{ j_1,\dots,j_m \}]^{\infty}) 
 \cdots
{{(-1)}^{k+1}}\bigcup P_n(C^A[\{ j_1,\dots,j_k \}]^{\infty}), 
\end{align*}
where
${{(-1)}^{m+1}} \bigcup_{\{j_1,\dots,j_m\} \subset J}$
means 
$\bigcup_{\{j_1,\dots,j_m\} \subset J}$ if $m$ is odd.
\end{lemma}
Hence we have
\begin{proposition}
Let $J =\{ j_1,\dots,j_k \}$ be a support subset of 
$\{1,2,\dots,N\}$. Then we have
\begin{align*}
&\zeta({C^A}^{\infty}, z)\\ 
= &\prod_{l=1}^k \zeta(C^A[\{ j_l \} ]^{\infty}, z)
\cdot
\prod_{\{j_1,j_2\} \subset J}
\zeta(C^A[\{ j_1,j_2 \}]^{\infty}, z)^{-1} \\
 \cdots
& \prod_{\{j_1,\dots,j_m\} \subset J}
\zeta(C^A[\{ j_1,\dots,j_m \}]^{\infty}, z)^{{(-1)}^{m+1}} 
 \cdots
\zeta(C^A[\{ j_1,\dots,j_k \}]^{\infty}, z)^{{(-1)}^{k+1}}. 
\end{align*}
\end{proposition}
\begin{corollary}
Suppose that there exists $j_0 \in \{1,2,\dots,N\}$ such that
$A(i,j_0) =1$ for all $i=1,2,\dots,N$.
Then 
$\zeta({C^A}^{\infty}, z) = f_{j_0}^A(z^2).
$
\end{corollary}
We reach the following formula of the zeta function of a Markov-Dyck shift of vertex type.
\begin{theorem}\label{thm:main}
Let $A$ be an $N \times N$ essential matrix with entries in $\{0,1\}$.
Then the zeta function $\zeta_{D_A}(z)$ of the Markov-Dyck shift $D_A$ is given by the following formula:
\begin{equation}
\zeta_{D_A}(z) = 
\frac{\zeta({C^A}^{\infty}, z)}{
\det(I_N - \diag[f_1^A(z^2),\dots,f_N^A(z^2)] zA)^2}
\end{equation}
where
\begin{equation*}
\zeta({C^A}^{\infty}, z) 
= \prod_{\{j_1,\dots,j_k\} \subset\{1,2,\dots,N\}}
   \zeta(C^A[\{ j_1,\dots,j_k \}]^{\infty}, z)^{{(-1)}^{k+1}}, 
\end{equation*}
the products
$\prod_{\{j_1,\dots,j_k\} \subset\{1,2,\dots,N\}}$
run over all subsets of $\{1,2,\dots,N\}$,
and
the zeta function
  $
 \zeta(C^A[\{ j_1,\dots,j_k \}]^{\infty}, z)
 $ is given by 
\begin{equation*}    
 \zeta(C^A[\{ j_1,\dots,j_k \}]^{\infty}, z) 
= 
\frac{1}{1-g_{C^A[\{ j_1,\dots,j_k \}]}(z)},
\end{equation*}
where
\begin{equation*}
g_{C^A[\{ j_1,\dots,j_k \}]}(z)
=  
z^2 \sum_{i=1}^N A(i,j_1)\cdots A(i,j_k) f_i^A(z^2),
\end{equation*}
and the functions 
$f_i^A(z^2), i=1,2,\dots, N$
satisfiy
the relations \eqref{eq:fia}.
\end{theorem}
\begin{proof}
For $n, k \in \N$, we define the following set 
$C_{n,k}^{A,-}$ similarly to $C_{n,k}^{A,+}$ by 
\begin{align*}
C_{n,k}^{A,-} = \{
(\alpha_{i_1},\dots,\alpha_{i_k},w_1,\dots,w_{2n})& \in B_{2n+k}(D_A)\mid \\
(w_1,\dots,w_{2n}) & \in C_n^A,\,
(\alpha_{i_1},\dots,\alpha_{i_k}) \in B_{k}(X_{A^t})
\}.
\end{align*} 
Similarly to the previous discussion,
we have a circular Markov code  
$
{\mathcal{C}}_A^-
=(C_A^-, I, \tilde{A^t}, r)
$ 
and the formula
\eqref{eq:zetacaplus}
for 
$
\zeta(\Omega_{{\mathcal{C}}_A^-},z).
$
We then have a disjoint union of periodic points
\begin{equation*}
P_n(D_A) 
=
P_n(\Omega_{{\mathcal{C}}_A^+}) \cup
P_n(\Omega_{{\mathcal{C}}_A^-}) \cup
P_n({C^A}^{\infty}) \cup
P_n(X_A) \cup
P_n(X_{A^t}).
\end{equation*}
Since 
$\zeta(\Omega_{{\mathcal{C}}_A^+},z)
=\zeta(\Omega_{{\mathcal{C}}_A^-},z),
$
Proposition \ref{prop:CAplus}
ensures us 
\begin{align*}
\zeta_{D_A}(z)
& = \zeta(\Omega_{{\mathcal{C}}_A^+},z)
\cdot
  \zeta(\Omega_{{\mathcal{C}}_A^-},z)
\cdot
  \zeta({C^A}^{\infty}, z)
\cdot
\frac{1}{\det(I_N - zA)}  
\cdot
\frac{1}{\det(I_N - zA^t)}\\
&=
\frac{\zeta({C^A}^{\infty}, z)}{
\det(I_N - \diag[f_1^A(z^2),\dots,f_N^A(z^2)] zA)^2}.
\end{align*}
\end{proof}
For a finite directed graph $G=(V,E)$ the above formula 
gives us the formula for the zeta function of the Markov-Dyck shift of 
vertex type.
\begin{corollary}\label{cor:zetaj0}
Suppose that there exists $j_0 \in \{1,2,\dots,N\}$ such that
$A(i,j_0) =$ for all $i=1,2,\dots,N$.
Then 
\begin{equation*}
\zeta_{D_A}(z)= 
\frac{f_{j_0}^A(z^2)}{
\det(I_N - \diag[f_1^A(z^2),\dots,f_N^A(z^2)] zA)^2}.
\end{equation*}
\end{corollary}
%%%%%%%%%%%%%%%%%%%%%%%%%%%%%%%%%%%%%%%%%%
\section{The zeta functions of Markov-Dyck shifts of edge type}
%%%%%%%%%%%%%%%%%%%%%%%%%%%%%%%%%%%%%%%%%%%%
The Markov-Dyck shifts in the paper 
\cite{KMMunster} are the 
Markov-Dyck shifts of edge type.
In \cite{KMMunster},
a formula of the zeta functions of Markov-Dyck shifts of edge type
has been presented.
In this section, we present the formula
\cite[Theorem 2.3]{KMMunster} from Theorem \ref{thm:main}.
We need the following lemma.
\begin{lemma}
For a finite directed graph $G =(V,E)$ with 
$|V| = N_0$ and $|E| = N_1$.
Let $f_1^V(x),\cdots,f_{N_0}^V(x)$
and $f_1^E(x),\cdots,f_{N_1}^E(x)$
be the functions satisfying the relations respectively
\begin{align} 
f_i^V(z) &
=  1 + z f_i^V(z) \sum_{j=1}^{N_0} A_G(j,i) f_j^V(z). \label{eq:fiv}\\
f_i^E(z) &
=  1 + z f_i^E(z) \sum_{j=1}^{N_1} A^G(j,i) f_j^E(z). \label{eq:fie}  
\end{align}
Then we have
\begin{align*}
& \det(I_{N_0} - \diag[f_1^V(z^2),\dots,f_{N_0}^V(z^2)] zA_G)\\
=
&\det(I_{N_1} - \diag[f_1^E(z^2),\dots,f_{N_1}^E(z^2)] zA^G).
\end{align*}
\end{lemma}
\begin{proof}
%%%%%%%%%%%%%%
\def\D^V{{ {D^V(z^2)} }}
\def\D^E{{ {D^E(z^2)} }}
%%%%%%%%%%%%%%%%%%%
Put the sets
$I_0 =\{1,2,\dots,N_0\}$,
$I_1 =\{1,2,\dots,N_1\}$
and
the diagonal matrices
$D^V(z^2) =\diag[f_1^V(z^2),\dots,f_{N_0}^V(z^2)]$
and
$D^E(z^2) =\diag[f_1^E(z^2),\dots,f_{N_1}^E(z^2)]$.
Define 
the $N_0 \times N_1$ matrix 
$S = [S(i,j)]_{i\in I_0, j \in I_1}$ and 
the $N_1 \times N_0$ matrix 
$R = [R(j,i)]_{j\in I_1, i \in I_0}$ by  
\begin{equation*}
S(i,j) =
\begin{cases}
1 & \text{ if } v_i = s(e_j),\\
0 & \text{ otherwise},
\end{cases}
\qquad
R(j,i) =
\begin{cases}
1 & \text{ if } t(e_j) = v_i,\\
0 & \text{ otherwise,}
\end{cases}
\end{equation*}
so that
$A_G = SR$ and $A^G = RS$.
For a vertex $v_i \in V$  
and en edge $e_j \in E$,
we set
\begin{align*}
C_n^{A_G}(v_i)
& = \{ (w_1,\dots,w_{2n}) \in C_n^{A_G} \mid 
       (v_i^*, w_1,\dots,w_{2n}, v_i) \in C_{n+1}^{A_G} \}, \\
C_n^{A^G}(e_j)
& = \{ (g_1,\dots,g_{2n}) \in C_n^{A^G} \mid 
       (e_j^*, g_1,\dots,g_{2n}, e_j) \in C_{n+1}^{A^G} \}.
\end{align*}
Let us denote by
$c_n^G(v_i)$ and 
$c_n^G(e_j)$
their cardinalities
$|C_n^{A_G}(v_i)|$
and
$|C_n^{A^G}(e_j)|$
respectively (\cite[pages 8,9]{MatsumotoIJM}).
Then we have
\begin{equation*}
f_i^V(z) = \sum_{n=0}^\infty c_n^G(v_i) z^n,
\qquad
f_j^E(z) = \sum_{n=0}^\infty c_n^G(e_j) z^n
\end{equation*} 
so that
$f_j^E(z) = f_i^V(z)$ when $s(e_j) = v_i$. 
Hence we have
\begin{equation*}
f_i^V(z^2) S(i,j) = S(i,j) f_j^E(z^2)
\end{equation*} 
which implies that
$D^V(z^2) S =S D^E(z^2).$
It then follows that
\begin{align*}
z D^V(z^2) A_G & = z D^V(z^2) S R = z S\cdot D^E(z^2) R, \\
z D^E(z^2) A^G & = z D^E(z^2) R S = D^E(z^2) R \cdot z S.
\end{align*} 
Hence 
the matrices $z D^V(z^2) A_G$ and $z D^E(z^2) A^G$ 
are elementary equivalent 
(see \cite[Definition 7.2.1]{LM}),
so that $\det(I_{N_0} - z D^V(z^2) A_G)
=
\det(I_{N_1} - z D^E(z^2) A^G).
$
\end{proof}
Therefore we have 
\begin{proposition}[{\cite[Theorem 2.3]{KMMunster}}]
If a matrix $A$ is an edge matrix $A^G=[A^G(e,f)]_{e,f \in E}$
defined by a finite directed graph $G =(V,E)$ with $|V| = N_0$,
then the zeta function of the Markov-Dyck shift 
$D_G (=D_{A^G})$ of edge type  
is given by the following formula:
\begin{equation}
\zeta_{D_G}(z) = 
\frac{\Pi_{i=1}^{N_0} f_i^G(z^2)}{
\det(I_N - \diag[f_1^G(z^2),\dots,f_{N_0}^G(z^2)] z A_G)^2} \label{eq:zetadg}
\end{equation}
where
$f_1^G(z^2),\dots,f_{N_0}^G(z^2)$ are the functions satisfying
\begin{equation}
f_i^G(z) 
=  1 + z f_i^G(z) \sum_{j=1}^{N_0} A_G(j,i) f_j^G(z). \label{eq:fig}  
\end{equation}
\end{proposition}
\begin{proof}
Since $f_i^G(x) = f_i^V(x), i=1,\dots,N_0$
and 
\begin{equation*}
\zeta({C^A}^{\infty}, z)
=\prod_{i=1}^{N_0}
\frac{1}{1 -g_{C^A[\{ j \}]}(z)}
=\prod_{i=1}^{N_0} f_i^G(z^2)
\end{equation*}
(cf. (cf. \cite[Proposition 4.7.11]{Stanley}),
the preceding lemma implies the equality
\eqref{eq:zetadg}
 from Theorem \ref{thm:main}.
\end{proof}
%%%%%%%%%%%%%%%%
%%%%%%%%%%%%%%%%%%%%%%%%%%%
\section{The Fibonacci-Dyck shift of vertex type}
%%%%%%%%%%%%%%%%%%%%%%%%%%

Let $G_2$ be the finite directed graph defined in the Figure 2.
%%%%%%%%%%%%%%%%%%%%%%%%%%%%%%%%%%%%%%%%%
\begin{figure}[htbp]
\begin{center}
%WinTpicVersion3.08
\unitlength 0.1in
\begin{picture}( 21.9100,  6.9000)( 16.5000,-24.0600)
% ELLIPSE 2 0 3 0
% 4 2092 2063 2170 2128 2170 2058 2170 2058
% 
\special{pn 8}%
\special{ar 2092 2064 78 66  0.0000000 6.2831853}%
% ELLIPSE 2 0 3 0
% 4 3762 2058 3841 2124 3841 2054 3841 2054
% 
\special{pn 8}%
\special{ar 3762 2058 80 66  0.0000000 6.2831853}%
% ELLIPSE 2 0 3 0
% 4 1832 2058 2014 2211 2003 2028 2056 2120
% 
\special{pn 8}%
\special{ar 1832 2058 182 154  0.3190696 6.0756891}%
% SARROW 2 0 3 1
% 2 2002 2112 2005 2106
% 
\special{pn 8}%
\special{pa 2002 2112}%
\special{pa 2006 2106}%
\special{fp}%
\special{sh 1}%
\special{pa 2006 2106}%
\special{pa 1958 2158}%
\special{pa 1982 2154}%
\special{pa 1994 2176}%
\special{pa 2006 2106}%
\special{fp}%
% ELLIPSE 2 0 3 0
% 4 2945 2128 2128 1850 2128 2124 3627 2120
% 
\special{pn 8}%
\special{ar 2946 2128 818 278  6.2480092 6.2831853}%
\special{ar 2946 2128 818 278  0.0000000 3.1562795}%
% SARROW 2 0 3 1
% 2 3762 2122 3762 2118
% 
\special{pn 8}%
\special{pa 3762 2122}%
\special{pa 3762 2118}%
\special{fp}%
\special{sh 1}%
\special{pa 3762 2118}%
\special{pa 3742 2186}%
\special{pa 3762 2172}%
\special{pa 3782 2186}%
\special{pa 3762 2118}%
\special{fp}%
% SPLINE 2 0 3 0
% 42 2133 2006 2133 1984 2139 1962 2150 1940 2166 1918 2187 1896 2212 1876 2242 1857 2276 1838 2316 1820 2358 1802 2406 1787 2455 1774 2509 1760 2565 1749 2625 1739 2686 1731 2748 1725 2812 1720 2877 1718 2941 1716 3006 1718 3071 1719 3136 1724 3199 1729 3260 1738 3319 1746 3376 1758 3430 1770 3482 1784 3529 1798 3574 1815 3614 1833 3650 1851 3680 1870 3708 1891 3729 1912 3747 1934 3759 1956 3766 1978 3767 2000 3767 2000
% 
\special{pn 8}%
\special{pa 2134 2006}%
\special{pa 2136 1974}%
\special{pa 2148 1946}%
\special{pa 2166 1918}%
\special{pa 2188 1896}%
\special{pa 2214 1876}%
\special{pa 2240 1858}%
\special{pa 2268 1842}%
\special{pa 2296 1830}%
\special{pa 2326 1816}%
\special{pa 2356 1804}%
\special{pa 2386 1794}%
\special{pa 2416 1784}%
\special{pa 2448 1776}%
\special{pa 2478 1768}%
\special{pa 2510 1760}%
\special{pa 2542 1754}%
\special{pa 2572 1748}%
\special{pa 2604 1742}%
\special{pa 2636 1738}%
\special{pa 2668 1734}%
\special{pa 2700 1730}%
\special{pa 2732 1728}%
\special{pa 2762 1724}%
\special{pa 2794 1722}%
\special{pa 2826 1720}%
\special{pa 2858 1720}%
\special{pa 2890 1718}%
\special{pa 2922 1716}%
\special{pa 2954 1716}%
\special{pa 2986 1718}%
\special{pa 3018 1718}%
\special{pa 3050 1718}%
\special{pa 3082 1720}%
\special{pa 3114 1722}%
\special{pa 3146 1726}%
\special{pa 3178 1728}%
\special{pa 3210 1730}%
\special{pa 3242 1736}%
\special{pa 3274 1740}%
\special{pa 3306 1744}%
\special{pa 3336 1750}%
\special{pa 3368 1756}%
\special{pa 3400 1764}%
\special{pa 3430 1770}%
\special{pa 3462 1778}%
\special{pa 3492 1788}%
\special{pa 3522 1796}%
\special{pa 3554 1806}%
\special{pa 3582 1820}%
\special{pa 3612 1832}%
\special{pa 3640 1846}%
\special{pa 3668 1862}%
\special{pa 3694 1880}%
\special{pa 3718 1902}%
\special{pa 3740 1926}%
\special{pa 3758 1952}%
\special{pa 3768 1982}%
\special{pa 3768 2000}%
\special{sp}%
% SARROW 2 0 3 1
% 2 2133 1992 2133 2006
% 
\special{pn 8}%
\special{pa 2134 1992}%
\special{pa 2134 2006}%
\special{fp}%
\special{sh 1}%
\special{pa 2134 2006}%
\special{pa 2154 1940}%
\special{pa 2134 1954}%
\special{pa 2114 1940}%
\special{pa 2134 2006}%
\special{fp}%
% STR 2 0 3 0
% 3 2090 1970 2090 2070 5 0
% 1
\put(20.9000,-20.7000){\makebox(0,0){1}}%
% STR 2 0 3 0
% 3 3760 1950 3760 2050 5 0
% 2
\put(37.6000,-20.5000){\makebox(0,0){2}}%
\end{picture}%
\end{center}
\caption{}
\end{figure}
%%%%%%%%%%%%%%%%%%%%%%%%%%%%%%%%%%%%%%%%%%%%%%%%%%%%%
The edge matrix $A^{G_2}$
and the vertex matrix $A_{G_2}$
are written as 
\begin{equation}
A^{G_2} =
\begin{bmatrix}
1 & 1 & 0 \\
0 & 0 & 1 \\
1 & 1 & 0   
\end{bmatrix},
\qquad
A_{G_2} =
\begin{bmatrix}
1 & 1 \\
1 & 0 
\end{bmatrix}
\end{equation}
respectively. 
We then have
\begin{proposition}
$D_{G_2}^V$ is not topologically conjugate to 
$D_{G_2}^E$.
\end{proposition}
\begin{proof}
It is easy to see that
the number of the $2$-periodic points of $D_{G_2}^V$ is $6$,
whereas 
that  of $D_{G_2}^E$ is $7$.
\end{proof}
The Fibonacci-Dyck shift $D^E_{G_2}$ 
of edge type is 
a subshift $D_{A^{G_2}}$
over six symbols which correspond to the edges of 
the directed graphs $G_2$ and $G_2^*$  of Figure 2.
The Fibonacci-Dyck shift $D_{A_{G_2}}$  
of vertex type is 
a subshift $D_{A_{G_2}}$
over four symbols which correspond to the vertices of 
the directed graphs of $G_2$ and $G_2^*$  of Figure 2.
Let us denote by 
$\alpha_1, \alpha_2$ and $\beta_1, \beta_2$
the symbols of $D_{A_{G_2}}$.
They have the following algebraic relations
from the relations \eqref{eq:CK} of operators
for
$A = A_{G_2}= 
\begin{bmatrix}
1 & 1 \\
1 & 0 
\end{bmatrix}:
$
\begin{equation*}
\alpha_1 \beta_1 =  \beta_1 \alpha_1 + \beta_2 \alpha_2=1,\qquad
\alpha_2 \beta_2 =  \beta_1 \alpha_1,\qquad
\beta_2 \alpha_2\beta_2 = \beta_2,
\end{equation*}
A word $\gamma=(\gamma_1,\dots,\gamma_m)$ of 
$\Sigma=\{\alpha_1, \alpha_2, \beta_1, \beta_2\}$
is forbidden if $\gamma_1\cdots\gamma_m=0$.
The Fibonacci-Dyck shift $D_{A_{G_2}}$ of vertex type
is defined as a subshift over $\Sigma$ whose forbidden words are defined 
in this sense.

We will compute the zeta function $\zeta_{D_{G_2}^V}(z)$
by using Corollary \ref{cor:zetaj0}.
Let $f_1(z), f_2(z)$ be the functions $f_1^V(z), f_2^V(z)$
which satisfy the following relations:
\begin{align*}
f_1(z) -1 & = z (f_1(z) + f_2(z)) f_1(z),\\
f_2(z) -1 & = z f_1(z)  f_2(z)   
\end{align*}
so that the equalities
\begin{equation*}
f_2(z)^2  = f_1(z),\qquad
z f_2(z)^3  -f_2(z) + 1  = 0    
\end{equation*}
hold (see \cite[Section 7]{MatsumotoIJM}).
We then have
\begin{align*}
\det(I_2 - \diag[f_1(z^2),f_2(z^2)]z A_{G_2})
& = \det\left(
{\begin{bmatrix}
1-zf_1(z^2) & -z f_1(z^2) \\
-z f_1(z^2) & 1 
\end{bmatrix}}\right)\\
& = 1 - zf_1(z^2)- z^2f_1(z^2)f_2(z^2) \\
& = 2 - zf_1(z^2)- f_2(z^2).  
\end{align*} 
\begin{proposition}
The zeta function $\zeta_{D_{G_2}^V}(z)$ 
of the Fibonacci-Dyck shift of vertex type
is 
\begin{equation}
\zeta_{D_{G_2}^V}(z) = \frac{1}{(2\xi(z)^2 + \xi(z) -1)^2}
\end{equation}  
where 
$\xi(z) = 
\frac{2}{\sqrt{3}} \sin(\frac{1}{3}\arcsin\frac{3\sqrt{3}}{2}z)
$ for
$
0 \le z \le \frac{2}{3\sqrt{3}}$.
\end{proposition}
\begin{proof}
By Corollary \ref{cor:zetaj0} with the above discussions, 
we have
\begin{align*}
\zeta_{D_{G_2}^V}(z)
& = \frac{f_1(z^2)}{(2 - z f_1(z^2) - f_2(z^2))^2}\\
& = \left( 
\frac{f_2(z^2)}{
(2f_2(z^2)-2 z^2 (f_2(z^2)^3)- z f_2(z^2)^2 - f_2(z^2))} \right)^2 \\
& = \frac{1}{
(1- 2 (z f_2(z^2)^2 - z f_2(z^2) )^2 }.
\end{align*}  
By putting $\xi(z) = z f_2(z^2)$,
we have
\begin{equation}
\zeta_{D_{G_2}^V}(z) = \frac{1}{(2\xi(z)^2 + \xi(z) -1)^2} \label{eq:zetaFibbV}
\end{equation}  
and
$\xi(z)^3 - \xi(z) +z =0$.
As in \cite[(4.10), (4.13)]{KMMunster},
we have
\begin{equation*}
\xi(z) = 
\frac{2}{\sqrt{3}} \sin(\frac{1}{3}\arcsin\frac{3\sqrt{3}}{2}z)
\quad
\text{ for }\quad
0 \le z \le \frac{2}{3\sqrt{3}}.
\end{equation*}
\end{proof}
We remark that 
the zeta function $\zeta_{D_{G_2}^E}(z)$ 
of the Fibonacci-Dyck shift of edge type
is 
\begin{equation*}
\zeta_{D_{G_2}^E}(z) = \frac{\xi(z)}{z (2\xi(z)^2 + \xi(z) -1)^2}
\qquad(\cite[Section 7]{KMMunster})
\end{equation*}  
which is different from \eqref{eq:zetaFibbV}.
%%%%%%%%%%%%%%%%%%%%%%%%%%%%%%%%%%

\medskip

{\it Acknowledgments:}
The author would like to thank 
Wolfgang Kriegerfor his various suggestions, comments, discussions 
and constant encouragements.
This work was supported by JSPS KAKENHI Grant Numbers 23540237.

\end{document}